\documentclass[a4paper,12pt]{amsart}
\usepackage[utf8]{inputenc}
\usepackage[margin=28mm,headheight=14pt]{geometry}
\usepackage{amsmath,amssymb,mathrsfs}
\usepackage{microtype,needspace}
\usepackage{tikz-cd}
\usepackage[hidelinks]{hyperref}

\numberwithin{equation}{section}
\allowdisplaybreaks[2]
\theoremstyle{plain}
\newtheorem{theorem}{Theorem}[section]
\newtheorem{proposition}[theorem]{Proposition}
\newtheorem{lemma}[theorem]{Lemma}
\newtheorem{corollary}[theorem]{Corollary}
\theoremstyle{definition}
\newtheorem{definition}[theorem]{Definition}
\theoremstyle{remark}
\newtheorem{remark}[theorem]{Remark}

\newcommand{\HC}{\mathrm{HC}}
\newcommand{\Lie}{\mathrm{Lie}}

\begin{document}
\title[Casselman's comparison conjecture]{A proof of Casselman's comparison conjecture}

\author{Junyan Xiao}
\address{(Xiao)School of Mathematical Sciences, Peking University, No. 5 Yiheyuan Road, Beijing 100871, P.R. China}
\email{jyxiao25@stu.pku.edu.cn}

\author{Jun Yu}
\address{(Yu)Jun Yu, School of Mathematical Sciences and Beijing International Center for Mathematical Research, Peking University, No. 5 Yiheyuan Road, Beijing 100871, China}
\email{junyu@bicmr.pku.edu.cn}

\begin{abstract}
We prove Casselman's comparison conjecture for parabolic subgroups of real linear reductive groups.
Let $Q=LU$ be a standard parabolic subgroup of a real linear reductive group $G$, with Levi subgroup $L$ and unipotent radical $U$.
For every Casselman--Wallach representation $\pi$ of $G$ and every $i\geq0$, the natural inclusion $\pi_K\hookrightarrow\pi$ induces a topological $L$-equivariant isomorphism
\[
 H_i(\mathfrak u,\pi_K)^\infty\xrightarrow{\ \sim\ }
 H_i(\mathfrak u,\pi).
\]
Put $\mathscr U=\mathcal U(\mathfrak u)$ and $F_k=\mathscr U/I_{\mathfrak u}^k$, where $I_{\mathfrak u}$ is the augmentation ideal.
For every $k\geq1$ and $i\geq0$, we establish natural topological $L$-isomorphisms
\[
\begin{aligned}
 \operatorname{Tor}_i^{\mathscr U}(F_k,\pi_K)^\infty
 \simeq \operatorname{Tor}_i^{\mathscr U}(F_k,\pi),\ 
 \operatorname{Ext}^i_{\mathscr U}(F_k,\pi_K)^\infty
 \simeq \operatorname{Ext}^i_{\mathscr U}(F_k,\pi).
\end{aligned}
\]
We also obtain closedness of $\mathfrak u^kM$
for the (co)chain spaces $M=C_j(\mathfrak u,\pi),C^j(\mathfrak u,\pi)$, and for $M=J_a,J_{a-1}/J_a$, where $(J_a)$ is the smooth Bruhat filtration of a principal series representation.
\end{abstract}
\maketitle
\setcounter{tocdepth}{2}
\begingroup
\small
\tableofcontents
\endgroup
\newpage

\section{Introduction}\label{S:introduction}

Let $\mathbf G$ be a connected linear reductive group over $\mathbb R$, and $G$ be a finite-index subgroup of $\mathbf G(\mathbb R)$ containing $\mathbf G(\mathbb R)^\circ$.
Fix a Cartan involution $\theta$ of $G$ and put $K=G^{\theta}$.

Let $P=MAN$ be a standard minimal parabolic subgroup of $G$.
That means, $A$ is a split torus and the Lie algebra of $A$ is a maximal abelian subspace of $\Lie G^{-\theta}$, $M=Z_{K}(A)$, and $\Lie N$ is the sum of positive root spaces from a chosen positive system $\Delta^{+}(\Lie G,\Lie A)$ of the root system $\Delta(\Lie G,\Lie A)$.

Let $Q=LU\supset P$ be a standard parabolic subgroup, with Levi subgroup $L\supset MA$ and unipotent radical $U\subset N$.
Then, $K_{L}:=K\cap L$ is a maximal compact subgroup of $L$.
Let the complexified Lie algebras be denoted by 
\[
 \mathfrak g=\operatorname{Lie}(G)\otimes_{\mathbb R}\mathbb C,
 \qquad \mathfrak l=\operatorname{Lie}(L)\otimes_{\mathbb R}\mathbb C,
 \qquad \mathfrak u=\operatorname{Lie}(U)\otimes_{\mathbb R}\mathbb C.
\]

For a complex Lie algebra $\mathfrak h$ and an $\mathfrak h$-module $V$, write $C_\bullet(\mathfrak h,V)=\bigwedge^\bullet\mathfrak h\otimes V$ for its Chevalley--Eilenberg chain complex and $H_i(\mathfrak h,V)$ for its homology.
Similarly, let  $C^\bullet(\mathfrak h,V)$ and $H^i(\mathfrak h,V)$ denote the cochain complex and Lie algebra cohomology.
For topological modules, homology and cohomology carry their quotient topologies.

For a Casselman--Wallach representation $\pi$ of $G$, let $\pi_K$ be its Harish--Chandra module of $K$-finite vectors.
The module $H_i(\mathfrak u,\pi_K)$ is Harish--Chandra for $(\mathfrak l,K_L)$; see \cite[Proposition~2.24]{HS83}.

For a Harish--Chandra $(\mathfrak l,K_L)$-module $W$, write $W^\infty$ for its Casselman--Wallach globalization, whose existence and uniqueness are guaranteed by the Casselman--Wallach globalization theorem \cite{Casselman89} and \cite[Chapter~11]{Wallach2}.
We prove the following statement.

\begin{theorem}\label{comparison conj}
For every $i\geq0$, the inclusion $\pi_K\hookrightarrow\pi$ induces a natural topological $L$-equivariant isomorphism
\[
 H_i(\mathfrak u,\pi_K)^\infty\xrightarrow{\ \sim\ }
 H_i(\mathfrak u,\pi).
\]
In particular, all differentials of $C_\bullet(\mathfrak u,\pi)$ have closed images, and its homology consists of Casselman--Wallach representations of $L$.
\end{theorem}

The cohomological version follows by Poincar\'e duality; its proof is given in Section~\ref{S:conclusion}.
\begin{corollary}\label{C:cohomological-comparison}
For every $i\geq0$, the inclusion $\pi_K\hookrightarrow\pi$ induces a natural topological $L$-equivariant isomorphism
\[
 H^i(\mathfrak u,\pi_K)^\infty\xrightarrow{\ \sim\ }
 H^i(\mathfrak u,\pi).
\]
All differentials of $C^\bullet(\mathfrak u,\pi)$ have closed images.
\end{corollary}

Casselman's comparison conjecture asks whether nilradical homology commutes with smooth globalization; Theorem~\ref{comparison conj} is its formulation for $Q$.
The broader comparison problem for canonical globalizations is discussed in \cite[Conjecture~10.3]{Vogan08}; see also the historical discussion in \cite[Introduction]{LLY}.

For a minimal parabolic, Hecht--Taylor \cite{Hecht-Taylor3} deduced the smooth comparison theorem from the corresponding result for minimal globalization.
Li--Liu--Yu \cite{LLY} gave a direct proof using distributions and the Bruhat and transverse-order filtrations, building on \cite{Casselman-Hecht-Milicic}.
The present proof combines comparison through Casselman--Jacquet modules with the distributional approach to treat arbitrary standard parabolics.
We also extend the comparison from trivial coefficients to higher augmentation quotients.

For a Lie algebra $\mathfrak h$, write $\mathcal U(\mathfrak h)$ for its universal enveloping algebra and $I_{\mathfrak h}=\ker(\mathcal U(\mathfrak h)\to\mathbb C)$ for its augmentation ideal.
For an $\mathfrak h$-module $M$, the notation $\mathfrak h^kM=I_{\mathfrak h}^kM$ denotes the algebraic span of products of $k$ Lie algebra operators, without taking closures.

Put $\mathscr U=\mathcal U(\mathfrak u),\
 F_k=\mathscr U/I_{\mathfrak u}^k\quad(k\geq1)$,
with the adjoint $L$-action on $F_k$.
For $k=1$, the following comparisons recover Theorem~\ref{comparison conj} and Corollary~\ref{C:cohomological-comparison}.

\begin{corollary}[Augmentation Tor and Ext comparison]\label{C:augmentation-comparison}\label{C:augmentation-Tor-comparison}\label{C:augmentation-Ext-comparison}
For every $k\geq1$ and $i\geq0$, the modules $\operatorname{Tor}_i^{\mathscr U}(F_k,\pi_K)$ and $\operatorname{Ext}^i_{\mathscr U}(F_k,\pi_K)$ are Harish--Chandra for $(\mathfrak l,K_L)$, and inclusion induces topological $L$-isomorphisms
\begin{equation}\label{Eq:augmentation-Tor-comparison}
 \operatorname{Tor}_i^{\mathscr U}(F_k,\pi_K)^\infty
 \xrightarrow{\ \sim\ }
 \operatorname{Tor}_i^{\mathscr U}(F_k,\pi),
\end{equation}
\begin{equation}\label{Eq:augmentation-Ext-comparison}
 \operatorname{Ext}^i_{\mathscr U}(F_k,\pi_K)^\infty
 \xrightarrow{\ \sim\ }
 \operatorname{Ext}^i_{\mathscr U}(F_k,\pi).
\end{equation}
The smooth modules carry the quotient topologies of the Chevalley--Eilenberg complexes described in Section~\ref{S:finite-coefficient-comparison}.
All differentials of these complexes have closed images, and their homology and cohomology are Casselman--Wallach representations of $L$.

The comparisons are natural in $\pi$ and compatible with the maps induced by $F_{k+1}\twoheadrightarrow F_k$, from level $k+1$ to $k$ for $\operatorname{Tor}$ and from level $k$ to $k+1$ for $\operatorname{Ext}$.
\end{corollary}

Finite-dimensional coefficients also give closedness in the chain and cochain spaces $C_j(\mathfrak u,\pi)=\bigwedge^j\mathfrak u\otimes\pi$ and $C^j(\mathfrak u,\pi)\simeq(\bigwedge^j\mathfrak u)^*\otimes\pi$.
We use the diagonal $\mathfrak u$-action, with the adjoint action on $\bigwedge^j\mathfrak u$ and its contragredient on $(\bigwedge^j\mathfrak u)^*$.

\begin{corollary}[Closedness in chain and cochain spaces]\label{C:Koszul-component-closedness}
For every $j\geq0$ and $k\geq1$, the algebraic subspaces
\[
 \mathfrak u^k C_j(\mathfrak u,\pi)
 \quad\text{and}\quad
 \mathfrak u^k C^j(\mathfrak u,\pi)
\]
are closed in $C_j(\mathfrak u,\pi)$ and $C^j(\mathfrak u,\pi)$, respectively.
\end{corollary}

We also obtain closedness for the smooth spaces arising from the Bruhat filtration.
Let $\sigma=\sigma_M\otimes e^\lambda$, where $\sigma_M$ is a finite-dimensional smooth representation of $M$ and $\lambda\in(\operatorname{Lie}(A)\otimes_{\mathbb R}\mathbb C)^*$ is arbitrary; let $N$ act trivially.
Put $X=G/P$, $E=G\times^P V_\sigma$, and $J=\operatorname{Ind}_P^G\sigma=C^\infty(X,E)$.
Here $V_\sigma$ is the representation space of $\sigma$.

Order the $Q$-orbits $O_1,\ldots,O_r$ in $X$ so that $Z_j=\bigcup_{a\leq j}O_a$ is closed, with $Z_0=\varnothing$, and put
\[
 J_j=\{f\in C^\infty(X,E):f\text{ is flat along }Z_j\},
 \qquad J_0=J,\quad J_r=0.
\]
The spaces $J_j$ carry their induced Fr\'echet topologies, and their successive quotients carry the quotient topologies.

\begin{corollary}[Closedness in smooth Bruhat pieces]\label{C:Bruhat-piece-closedness}
Let $M=J_j$ for $0\leq j\leq r$, or $M=J_{j-1}/J_j$ for $1\leq j\leq r$.
All differentials of $C_\bullet(\mathfrak u,M)$ have closed images, and its homology consists of Casselman--Wallach representations of $L$.
For every $k\geq1$, the algebraic subspace $\mathfrak u^kM$ is closed in $M$.
\end{corollary}

\subsection*{Conventions}
For a nilpotent Lie algebra $\mathfrak v$ and a module $M$, put
\[
 M^{[\mathfrak v]}=\bigcup_{k\geq1}
 \{m\in M:X_1\cdots X_km=0
       \text{ for all }X_1,\ldots,X_k\in\mathfrak v\}.
\]
This is the generalized-invariant submodule, corresponding to the zero character.
For an abelian Lie algebra $\mathfrak z$, we write $M^{\mathfrak z\text{-fin}}$ for the subspace of vectors generating finite-dimensional $\mathcal U(\mathfrak z)$-modules; it is the direct sum of all generalized $\mathfrak z$-weight spaces.
We use $E'$ for the continuous strong dual and $V^*$ for the algebraic dual.
For an affine real algebraic variety $Y$, $R(Y)$ is its complexified ring of regular functions.

A short exact sequence $0\to E_1\to E_2\to E_3\to0$ of Fr\'echet spaces is \emph{strict} if the first map is a homeomorphism onto a closed subspace and the second identifies $E_3$ with the topological quotient.
A functor is \emph{strictly exact} if it sends strict short exact sequences to strict short exact sequences.

\subsection*{Proof outline}
For a Harish--Chandra $(\mathfrak g,K)$-module $V$ and a Casselman--Wallach representation $\pi$ of $G$, the \emph{Casselman--Jacquet functor} and the \emph{smooth Casselman--Jacquet functor} are given by
\[
 \mathrm{CJ}_U(V)=\varprojlim_{k\geq1}V/\mathfrak u^kV,
 \qquad
 \mathrm{CJ}^{\infty}_U(\pi)
   =\varprojlim_{k\geq1}\pi/\overline{\mathfrak u^k\pi}.
\]
The bar denotes closure, and the smooth completion carries the inverse-limit topology.
Set $V=\pi_K$ and $E=\mathrm{CJ}_U(V)$.
Write $\widehat E$ for the globalization of the complete module $E$ constructed in Section~\ref{S:complete-homology}, and $q_\pi:\pi\to\mathrm{CJ}^{\infty}_U(\pi)$ for the canonical surjection \cite[Proposition~6.1]{CWYZ}.
The comparison is obtained from the following chain of isomorphisms:
\begin{equation}\label{Eq:intro-chain}
\begin{aligned}
 &H_i(\mathfrak u,V)^\infty\\
 &\xrightarrow{\sim} H_i(\mathfrak u,E)^\infty\\
 &\xrightarrow{\sim} H_i(\mathfrak u,\widehat E)\\
 &\xrightarrow{\sim} H_i(\mathfrak u,\mathrm{CJ}^{\infty}_U(\pi))\\
 &\xleftarrow[\sim]{q_{\pi,*}}
 H_i(\mathfrak u,\pi).
\end{aligned}
\end{equation}
The starting point is degree zero: passing from $V$ to its $\mathfrak u$-adic completion $E$ leaves the coinvariants $V/\mathfrak u V$ unchanged.
This suggests completing the complex that computes the higher homology.
Its completed chains are the chains for $E$, and exactness moves completion to homology, where $\mathfrak u$ already acts trivially (Section~\ref{S:algebraic-comparison}).

For the second comparison, ordinary globalization is available on each finite quotient of $E$.
We therefore look for pieces of $E$ that stabilize at a finite level.
A fixed central weight would have such property because $\mathfrak u$ raises the weight height.
This leads us to globalize the resulting weight subcomplex (Section~\ref{S:complete-homology}).
The third follows from the finite-level globalization theorem of \cite{CWYZ}.
These three comparisons are also shown in \cite{Wu-Yu}. 

It remains to prove that $q_{\pi,*}$ is an isomorphism. For a principal series $J$, duality reduces the assertion for $q_{J,*}$ to showing that inclusion of generalized invariants induces an isomorphism
\begin{equation}\label{Eq:intro-dual}
 H^i(\mathfrak u,(J')^{[\mathfrak u]})
 \xrightarrow{\ \sim\ } H^i(\mathfrak u,J'),
\end{equation}
which naturally leads to a distributional approach.
The Bruhat and normal-order filtrations reduce this comparison to the normal graded pieces on a $Q$-orbit (Section~\ref{S:filtrations}), which is further reduced to the computation of cohomology along unipotent fibers.
We prove this comparison by a variant of the central-subgroup induction of \cite[Section~4]{LLY} in Section~\ref{S:comparison-poly}.
Finally, general $\pi$ follows from coresolution by principal series (Theorem~\ref{T:smooth-CJ-comparison}).

In Section~\ref{S:finite-coefficients}, we extend comparison theorem to finite-dimensional coefficients.
Applying this to augmentation quotients gives Corollary~\ref{C:augmentation-comparison}; their zero-degree homology also detects the closedness in Corollary~\ref{C:Koszul-component-closedness}.

\subsection*{Acknowledgement}
Theorem \ref{comparison conj} and Corollary \ref{C:cohomological-comparison} are folklore theorems of Casselman. As presented in this paper, the proof of them consists of two main parts: one is the study of Casselman-Jacquent functors (\cite{CWYZ}, \cite{Wu-Yu} and Section \ref{S:CJ-homology} of this paper), another is the extension of the study of Bruhat filtration and the proof of the cohomological comparison between a Casselman-Wallach representations and its generalized $\mathfrak{u}$-invariant submodules. We had these ideas in the end of 2025 and wrote a draft in June 2026 without using AI. Jun Yu learned Corollaries 
\ref{C:augmentation-comparison}, \ref{C:Koszul-component-closedness}, \ref{C:Bruhat-piece-closedness} from Casselman during email communication in January 2021. They are shown here with AI-assistance. AI assistance was also used to audit the manuscript and revise proofs and exposition. Hang Xue notified us that he obtained different proofs of Theorem \ref{comparison conj}, 
Corollary \ref{C:cohomological-comparison} and 
Corollary \ref{C:augmentation-Tor-comparison} with AI-assistance (\cite{Xue}). Jun Yu thanks William Casselman, Gang Liu, Ning Li, Yoshiki Oshima, Binyong Sun, David Vogan and Hang Xue for helpful communications, and thanks Kei Yuen Chan, Kaidi Wu, Hongfeng Zhang for collaborations leading to \cite{CWYZ} and \cite{Wu-Yu}. Our research is supported by NSFC Grants 12431001 and 11971036.

\section{Completion and globalization}\label{S:CJ-homology}

Fix the standard parabolic $Q=LU$.
Put $\mathfrak q=\mathfrak l\ltimes\mathfrak u$ and let $\mathfrak z_L$ be the center of $\mathfrak l$.
We write $\HC(\mathfrak g,K)$ and $\HC(\mathfrak l,K_L)$ for the usual Harish--Chandra categories, $\mathrm{CW}_G$ and $\mathrm{CW}_L$ for the corresponding categories of Casselman--Wallach representations, and $W^\infty$ for the globalization of a Harish--Chandra module $W$.
The notation $\HC(\mathfrak g,K_L)$ refers to the complete modules defined below.

This section is an update of some parts of \cite{CWYZ} and \cite{Wu-Yu}. Several results of this section have already appeared in \cite{CWYZ} and \cite[Sections~4--5]{Wu-Yu}.
The complete category, its generalized central-weight decomposition, and its inverse-limit globalization are those of \cite[Definition~4.1, Lemma~4.4, Definition~4.16, and~(4.4)]{Wu-Yu}.
The algebraic comparison \eqref{Eq:algebraic-comparison} and the compatibility of globalization with the Casselman--Jacquet functor in Corollary~\ref{C:CJ-comparison} recover \cite[Propositions~5.2 and~4.18]{Wu-Yu}, respectively.
The homology comparison in Theorem~\ref{T:complete-global-comparison} also follows from \cite[Proposition~4.13 and Theorem~5.3]{Wu-Yu} by finite-length d\'evissage.
Here we give direct proofs using Artin--Rees completion and total-weight subcomplexes, and record the closed image and closed subquotient statements needed below.

\subsection{Homology and algebraic completion}\label{S:algebraic-comparison}

\begin{definition}\label{D:complete}
Let $\HC(\mathfrak g,K_L)$ be the category of vector spaces $E$ with compatible $\mathcal U(\mathfrak g)$- and $K_L$-actions such that
\begin{enumerate}
\item $E/\mathfrak uE$ is a Harish--Chandra $(\mathfrak l,K_L)$-module;
\item the action is $\mathfrak u$-complete: the canonical map
      \[
       E\longrightarrow\varprojlim_k E/\mathfrak u^kE
      \]
      is an isomorphism.
\end{enumerate}
Morphisms are linear maps commuting with both actions.
\end{definition}

Put $A=\mathcal U(\mathfrak q)$ and $I=\mathfrak uA=A\mathfrak u$, so $I^kM=\mathfrak u^kM$ for every $A$-module $M$.

\begin{lemma}\label{L:adic-completion}
The $I$-adic completion functor is exact on finitely generated $A$-modules.
Put $\widehat A=\varprojlim_k A/I^k$.
For every finitely generated $A$-module $M$, there is a natural isomorphism
\[
 M^{\wedge I}:=\varprojlim_k M/I^kM\simeq\widehat A\otimes_A M.
\]
Moreover, for every $k\geq1$, the completion map induces a natural isomorphism
\[
 M/I^kM\xrightarrow{\ \sim\ }M^{\wedge I}/I^kM^{\wedge I}.
\]
\end{lemma}

\begin{proof}
Since $\mathfrak u$ is a nilpotent ideal of $\mathfrak q$, $I$ has the strong Artin--Rees property \cite[Theorem~2.1]{SW82}.
Since $A$ is Noetherian, exactness of completion and the tensor-product description follow from \cite[Chapter~I, Proposition~5.10(1),(3)]{NVO79}.
Applying these assertions to $0\to I^k\to A\to A/I^k\to0$, and likewise to right $A$-modules, gives
\[
 \widehat A I^k=\ker(\widehat A\to A/I^k)=I^k\widehat A.
\]
It follows that
\[
 M^{\wedge I}/I^kM^{\wedge I}\simeq (\widehat A/I^k\widehat A)\otimes_A M\simeq M/I^kM .
\]
These identifications are natural and identify the map induced by completion with the asserted isomorphism.
\end{proof}

\begin{proposition}\label{P:completion-homology}
For every finitely generated $A$-module $M$ and every $i\geq0$, $H_i(\mathfrak u,M)$ is finitely generated over $\mathcal U(\mathfrak l)$, and the completion map induces a natural $\mathfrak l$-isomorphism
\[
 H_i(\mathfrak u,M)\xrightarrow{\ \sim\ } H_i(\mathfrak u,M^{\wedge I}).
\]
This isomorphism is $K_L$-equivariant whenever $M$ carries a compatible $K_L$-action.
\end{proposition}

\begin{proof}
Give $F\otimes T$ the diagonal action $\mathcal L$ for a finite-dimensional $\mathfrak q$-module $F$ and any $A$-module $T$.
For $X\in\mathfrak q$ one has
\[
 f\otimes Xv=\mathcal L_X(f\otimes v)-(Xf)\otimes v .
\]
If $\mathfrak u^{c+1}F=0$, then for all $k\geq1$,
\begin{equation}\label{Eq:tensor-cofinal}
 I^{k+c}(F\otimes T)\subset F\otimes\mathfrak u^kT,
 \qquad F\otimes\mathfrak u^{k+c}T\subset I^k(F\otimes T).
\end{equation}

Take $F_j=\bigwedge^j\mathfrak u$ and $C_j=F_j\otimes M$.
Then $C_j$, hence the homology is finitely generated over $A$.
For $X\in\mathfrak u$, define $s_X(f\otimes v)=X\wedge f\otimes v$, which satisfies $ds_X+s_Xd=\mathcal L_X$.
Thus $I$ annihilates homology, giving that $H_i(\mathfrak u,M)$ is finitely generated over $A/I=\mathcal U(\mathfrak l)$.

The adjoint $\mathfrak u$-action on $F_j$ is nilpotent.
By \eqref{Eq:tensor-cofinal}, the completion of $C_j$ identifies with $F_j\otimes M^{\wedge I}$.
Exactness of completion, by Lemma~\ref{L:adic-completion}, now gives
\begin{equation}\label{Eq:completed-CE-homology}
 H_i(\mathfrak u,M^{\wedge I})\simeq H_i(\mathfrak u,M)^{\wedge I}=H_i(\mathfrak u,M).
\end{equation}
Under these identifications the homology map induced by $M\to M^{\wedge I}$ is the canonical map to the completion of $H_i(\mathfrak u,M)$, and preserves $K_L$-action.
\end{proof}

For any $\mathfrak g$-module $T$, evaluation $\mathfrak g\otimes T\to T$ is $\mathfrak q$-equivariant.
The adjoint action of $\mathfrak u$ on $\mathfrak g$ is unipotent, so $\mathfrak u^{c+1}\mathfrak g=0$ for some $c$. By these and (the second inclusion of) \eqref{Eq:tensor-cofinal} we have
\begin{equation}\label{Eq:cross-level-bound}
 X\mathfrak u^{k+c}T\subset\mathfrak u^kT
 \qquad(X\in\mathfrak g,\ k\geq1).
\end{equation}

Now we apply Proposition~\ref{P:completion-homology} to $V\in\HC(\mathfrak g,K)$.
By Casselman--Osborne \cite{Casselman-Osborne}, \cite[Lemma~2.22]{HS83}, $V$ is finitely generated over $\mathcal U(\mathfrak n)$ and hence over $A$.
For $E=\mathrm{CJ}_U(V)=V^{\wedge I}$, Lemma~\ref{L:adic-completion} identifies the finite quotients as
\begin{equation}\label{Eq:CJ-algebraic-levels}
 E/\mathfrak u^kE\simeq V/\mathfrak u^kV\qquad(k\geq1).
\end{equation}
Thus $E$ is $\mathfrak u$-complete.
The $\mathfrak g$-action extends to $E$ by \eqref{Eq:cross-level-bound}, compatibly with the $K_L$-action.
Hecht--Schmid \cite[Proposition~2.24]{HS83} gives $H_i(\mathfrak u,V)\in\HC(\mathfrak l,K_L)$.
In degree zero this identifies $E/\mathfrak uE$ with the Harish--Chandra module $H_0(\mathfrak u,V)$, so $E$ lies in $\HC(\mathfrak g,K_L)$.
Proposition~\ref{P:completion-homology} gives the first comparison in $\HC(\mathfrak l,K_L)$,
\begin{equation}\label{Eq:algebraic-comparison}
 H_i(\mathfrak u,V)\xrightarrow{\ \sim\ }H_i(\mathfrak u,\mathrm{CJ}_U(V)).
\end{equation}

\subsection{Globalization by total weights}\label{S:complete-homology}
We recall the weight notation.
For a $\mathfrak z_L$-module $E$ and $\alpha\in\mathfrak z_L^*$, put
\[
 E_\alpha=\{v\in E:\text{ for some }r\geq1,
       (Z-\alpha(Z))^rv=0\text{ for all }Z\in\mathfrak z_L\}.
\]
Write $\mathrm{wt}(E)$ for the set of weights with nonzero weight space.
These weight spaces are $K_L$-stable since $L$ centralizes $\mathfrak z_L$.

\begin{lemma}\label{L:complete-weights}
Let $E\in\HC(\mathfrak g,K_L)$ and $M_k=E/\mathfrak u^kE$.
Each $M_k$ is Harish--Chandra for $(\mathfrak l,K_L)$.
For fixed $\alpha$, the projection identifies $E_\alpha$ with $(M_k)_\alpha$ for all sufficiently large $k$.
In particular, $E_\alpha$ is Harish--Chandra.
The weight projections induce an isomorphism commuting with the actions of $\mathfrak l$ and $K_L$,
\begin{equation}\label{Eq:complete-weight-product}
 E\simeq\prod_\alpha E_\alpha .
\end{equation}
The submodule $E_{\mathrm{fin}}=E^{\mathfrak z_L\text{-fin}}=\bigoplus_\alpha E_\alpha$ is finitely generated over $\mathcal U(\mathfrak q)$.
\end{lemma}
\begin{proof}
Multiplication realizes each graded piece $\mathfrak u^jE/\mathfrak u^{j+1}E$ as a quotient of $\mathfrak u^{\otimes j}\otimes M_1$.
Since $M_1$ is Harish--Chandra, so are the graded pieces and hence each finite extension $M_k$.

Choose $H_Q\in\mathfrak z_L$ with positive integral adjoint eigenvalues on $\mathfrak u$ and put $h(\alpha)=\operatorname{Re}\alpha(H_Q)$.
Let $c$ be a lower bound for the heights of the finite set $S=\mathrm{wt}(M_1)$.
Every weight of $\mathfrak u^jE/\mathfrak u^{j+1}E$ has height at least $c+j$, so $(M_{k+1})_\alpha\to(M_k)_\alpha$ is an isomorphism for $k>h(\alpha)-c$.

A vector at a stable level $(M_k)_\alpha$ lifts uniquely through $(M_l)_\alpha\ (l\geq k)$ by these isomorphism, so by completeness the lift belongs to $E_\alpha$.
Thus $E_\alpha\simeq \varprojlim_k (M_k)_\alpha$ is the stable weight space; for fixed $\alpha$ this is identified with $(M_k)_\alpha$ for $k$ sufficiently large, hence Harish--Chandra.
Taking inverse limits of the finite weight decompositions of the $M_k$ gives \eqref{Eq:complete-weight-product}.

It remains to prove finite generation.
If $\alpha\notin S$, then $E_\alpha\subset\mathfrak uE$.
Choose a weight basis $(X_a)$ of $\mathfrak u$, with weights $\gamma_a$.
For $v\in E_\alpha$, write $v=\sum_a X_ae_a$.
Projecting to weight $\alpha$ gives
\[
 v=\sum_a X_a\operatorname{pr}_{\alpha-\gamma_a}(e_a),\quad \operatorname{pr}_{\alpha-\gamma_a}:E\to E_{\alpha-\gamma_a}.
\]
Note that $\alpha-\gamma_a$ has height at most $h(\alpha)-1$.
Since all weights of $E$ have height at least $c$, a finite iteration shows that $\bigoplus_{\lambda\in S}E_\lambda$ generates $E_{\mathrm{fin}}$ over $\mathcal U(\mathfrak u)$.
This finite sum is Harish--Chandra and hence finitely generated over $\mathcal U(\mathfrak l)$, proving the assertion.
\end{proof}

For $E\in\HC(\mathfrak g,K_L)$, put $M_k=E/\mathfrak u^kE$ and define $\widehat E=\varprojlim_k M_k^\infty$ with its inverse-limit topology.
Each $M_k$ has only finitely many weights.
Globalizing its weight decomposition gives
\[
 M_k^\infty
 \simeq\bigoplus_{\alpha\in\mathrm{wt}(M_k)}((M_k)_\alpha)^\infty,
\]
and the transition maps preserve the weight summands.

Let $p_{k,\alpha}:E_\alpha\to(M_k)_\alpha$ be the natural projection.
A family $(v_\alpha)_\alpha\in\prod_\alpha(E_\alpha)^\infty$ determines an element of $\widehat E$ with $k$-th coordinate
\[
 v_k=\sum_{\alpha\in\mathrm{wt}(M_k)}p_{k,\alpha}^\infty(v_\alpha).
\]
These finite sums are compatible with the transition maps.

Conversely, let $(v_k)_k\in\widehat E$.
For each $\alpha$, Lemma~\ref{L:complete-weights} identifies $(M_k)_\alpha$ with $E_\alpha$ for all sufficiently large $k$.
The $\alpha$-component of $v_k$ at any such level therefore determines $v_\alpha\in(E_\alpha)^\infty$, independently of $k$.
Compatibility recovers every $v_k$ by the sum above.

Both maps are continuous.
Each $v_k$ depends continuously on finitely many weight coordinates, and each $v_\alpha$ is read off continuously at a single stable level.
We thus obtain a natural topological $L$-isomorphism
\begin{equation}\label{glob-iso}
 \widehat E\simeq\prod_\alpha(E_\alpha)^\infty.
\end{equation}
Under \eqref{Eq:complete-weight-product}, the map $E\to\widehat E$ is the product of the embeddings $E_\alpha\to(E_\alpha)^\infty$, and has dense image.

We next extend the $\mathfrak g$-action.
Choose $c$ as in \eqref{Eq:cross-level-bound}.
Since $\mathfrak g\,\mathfrak u^{k+c}E\subset\mathfrak u^kE$, the action on $E$ induces well-defined maps
\[
 a_k:\mathfrak g\otimes M_{k+c}\longrightarrow M_k,\qquad
 X\otimes[e]_{k+c}\longmapsto[Xe]_k.
\]
With the adjoint action on $\mathfrak g$ and the diagonal action on the tensor product, $a_k$ is $(\mathfrak l,K_L)$-equivariant.
The maps $a_k$ commute with the transition maps.
Hence $a_k$ globalizes to a continuous $L$-equivariant map $a_k^\infty:\mathfrak g\otimes M_{k+c}^\infty\to M_k^\infty$.
For $v=(v_k)_k\in\widehat E$, set
\[
 (Xv)_k=a_k^\infty(X\otimes v_{k+c}).
\]
The coordinates on the right are compatible, so this formula defines a continuous $L$-equivariant map $\mathfrak g\otimes\widehat E\to\widehat E$ extending the action on $E$.
This gives $\widehat E$ a compatible $(\mathfrak g,L)$-action.

\begin{theorem}\label{T:complete-global-comparison}
For every $E\in\HC(\mathfrak g,K_L)$ and $i\geq0$, the module $H_i(\mathfrak u,E)$ is Harish--Chandra for $(\mathfrak l,K_L)$, and $E\to\widehat E$ induces a natural topological $L$-isomorphism
\begin{equation}\label{Eq:complete-comparison-summary}
 H_i(\mathfrak u,E)^\infty\xrightarrow{\ \sim\ }H_i(\mathfrak u,\widehat E).
\end{equation}
Every differential of $C_\bullet(\mathfrak u,\widehat E)$ has closed image.
\end{theorem}
\begin{proof}
Put $F_j=\bigwedge^j\mathfrak u$.
The chains of total weight $\alpha$ form the space
\[
 C_{j,\alpha}=\bigoplus_{\beta\in\mathrm{wt}(F_j)}
        (F_j)_\beta\otimes E_{\alpha-\beta}.
\]
This is a finite sum of Harish--Chandra modules.
Indeed, each $E_{\alpha-\beta}$ is Harish--Chandra by Lemma~\ref{L:complete-weights}, and $(F_j)_\beta$ is a finite-dimensional $L$-module.

A direct computation shows that the Chevalley--Eilenberg differential preserves total weight.
Thus $C_{\bullet,\alpha}$ is a complex in $\HC(\mathfrak l,K_L)$, and
\[
 C_\bullet(\mathfrak u,E_{\mathrm{fin}})=\bigoplus_\alpha C_{\bullet,\alpha},
 \qquad C_\bullet(\mathfrak u,E)=\prod_\alpha C_{\bullet,\alpha}.
\]

Fix $i$ and put $D_{i,\alpha}=H_i(C_{\bullet,\alpha})$.
Taking homology in the direct-sum decomposition of $C_\bullet(\mathfrak u,E_{\mathrm{fin}})$ gives
\[
 H_i(\mathfrak u,E_{\mathrm{fin}})=\bigoplus_\alpha D_{i,\alpha}.
\]
By Lemma~\ref{L:complete-weights}, $E_{\mathrm{fin}}$ is finitely generated over $\mathcal U(\mathfrak q)$.
Proposition~\ref{P:completion-homology} therefore shows that this direct sum is finitely generated over $\mathcal U(\mathfrak l)$.
Since $\mathfrak l$ preserves each summand, $\Sigma_i=\{\alpha:D_{i,\alpha}\ne0\}$ is finite.
Taking homology in the product decomposition of $C_\bullet(\mathfrak u,E)$ now gives
\[
 H_i(\mathfrak u,E)=\prod_\alpha D_{i,\alpha}
                   =\bigoplus_{\alpha\in\Sigma_i}D_{i,\alpha}.
\]
Each $D_{i,\alpha}$ is Harish--Chandra, so the same holds for this finite sum.

By \eqref{glob-iso}, $C_\bullet(\mathfrak u,\widehat E)$ is the topological product of the complexes $(C_{\bullet,\alpha})^\infty$, and the differentials agree on the dense algebraic weight spaces.
Exactness of globalization on $\HC(\mathfrak l,K_L)$ gives $H_i((C_{\bullet,\alpha})^\infty)\simeq(D_{i,\alpha})^\infty$ topologically, with closed differential images.
The image of differential is a product of closed subspaces and is therefore closed.
Since $\Sigma_i$ is finite, we obtain
\[
 H_i(\mathfrak u,E)^\infty
 \simeq\left(\bigoplus_{\alpha\in\Sigma_i}D_{i,\alpha}\right)^\infty
 \simeq\bigoplus_{\alpha\in\Sigma_i}(D_{i,\alpha})^\infty
 \simeq H_i(\mathfrak u,\widehat E).
\]
On each summand the map induced by $E\to\widehat E$ is the canonical embedding $D_{i,\alpha}\to(D_{i,\alpha})^\infty$, so this is the isomorphism in \eqref{Eq:complete-comparison-summary}.
\end{proof}

\begin{remark}\label{R:weightwise-globalization}
The proof applies more generally to a finitely generated $(\mathfrak q,K_L)$-module $M=\bigoplus_\alpha M_\alpha$ whose generalized $\mathfrak z_L$-weight spaces $M_\alpha$ are Harish--Chandra for $(\mathfrak l,K_L)$.

Suppose that $\mathcal M=\prod_\alpha M_\alpha^\infty$, with the product topology and diagonal $L$-action, carries a compatible continuous $(\mathfrak q,L)$-action extending the $\mathfrak q$-action on $M$.

Then, for every $i\geq0$, $H_i(\mathfrak u,M)$ is Harish--Chandra, and inclusion induces a natural topological $L$-isomorphism
\[
 H_i(\mathfrak u,M)^\infty\xrightarrow{\ \sim\ }H_i(\mathfrak u,\mathcal M).
\]
The differentials of $C_\bullet(\mathfrak u,\mathcal M)$ have closed images, and its homology consists of Casselman--Wallach representations of $L$.

Indeed, the total-weight chain spaces are Harish--Chandra.
Proposition~\ref{P:completion-homology} makes $H_i(\mathfrak u,M)$ finitely generated over $\mathcal U(\mathfrak l)$, so only finitely many total weights contribute in each degree.
Exactness of globalization on each weight block and the product topology then give the assertions, as above.
\end{remark}

\begin{corollary}\label{C:CJ-comparison}
For $\pi\in\mathrm{CW}_G$ and $V=\pi_K$, the maps induced by $V\hookrightarrow\pi$ give a natural topological $(\mathfrak g,L)$-isomorphism
\[
 \widehat{\mathrm{CJ}_U(V)}\xrightarrow{\sim}\mathrm{CJ}^{\infty}_U(\pi).
\]
For every $i\geq0$, the composite $V\hookrightarrow\pi\xrightarrow{q_\pi}\mathrm{CJ}^{\infty}_U(\pi)$ induces a natural topological $L$-isomorphism
\begin{equation}\label{Eq:CJ-comparison}
 H_i(\mathfrak u,V)^\infty\xrightarrow{\ \sim\ }
 H_i(\mathfrak u,\mathrm{CJ}^{\infty}_U(\pi)).
\end{equation}
In particular, the differentials of $C_\bullet(\mathfrak u,\mathrm{CJ}^{\infty}_U(\pi))$ have closed images.
\end{corollary}
\begin{proof}
Put $E=\mathrm{CJ}_U(V)$.
By \eqref{Eq:CJ-algebraic-levels} and the finite-level globalization theorem \cite[Theorem~13.1 and Remark~13.2]{CWYZ}, the canonical maps give
\[
 (E/\mathfrak u^kE)^\infty\simeq(V/\mathfrak u^k V)^\infty
 \xrightarrow{\ \sim\ }\pi/\overline{\mathfrak u^k\pi},
\]
Taking inverse limits gives a topological $L$-isomorphism $\widehat{\mathrm{CJ}_U(V)}\simeq\mathrm{CJ}^{\infty}_U(\pi)$.
By \eqref{Eq:cross-level-bound} and continuity on $\pi$, each $X\in\mathfrak g$ acts continuously, and the isomorphism intertwines these actions on the dense image of $V$, hence everywhere.

Combined with \eqref{Eq:algebraic-comparison} after globalization and with Theorem~\ref{T:complete-global-comparison} gives \eqref{Eq:CJ-comparison}.
At each level the map from $V$ is induced by $V\hookrightarrow\pi$, proving the asserted compatibility.
\end{proof}

The weight decomposition also passes to closed subquotients, as follows.

\begin{lemma}[Closed subquotients]\label{L:closed-weight-subquotients}
Let $E\in\HC(\mathfrak g,K_L)$, and let $W'\subset W\subset\widehat E$ be closed $(\mathfrak q,L)$-submodules.
Put
\[
 W_{\mathrm{fin}}=W\cap E_{\mathrm{fin}},
 \qquad
 W'_{\mathrm{fin}}=W'\cap E_{\mathrm{fin}}.
\]
Then $W_{\mathrm{fin}}$ and $W_{\mathrm{fin}}/W'_{\mathrm{fin}}$ are finitely generated over $\mathcal U(\mathfrak q)$ and have Harish--Chandra generalized $\mathfrak z_L$-weight spaces.
There are natural topological $L$-isomorphisms
\[
 \begin{aligned}
 W&\simeq\prod_\alpha\bigl((W_{\mathrm{fin}})_\alpha\bigr)^\infty,\\
 W/W'&\simeq\prod_\alpha
       \bigl((W_{\mathrm{fin}}/W'_{\mathrm{fin}})_\alpha\bigr)^\infty.
 \end{aligned}
\]
\end{lemma}
\begin{proof}
Since $W$ is closed and $L$-stable, it is stable under $\mathcal U(\mathfrak z_L)$.
Finite generation of each $E_\beta$ gives a uniform generalized weight exponent, and the corresponding annihilating relations extend to $E_\beta^\infty$ by density.

For each finite set of weights $S\ni\alpha$, the Chinese remainder theorem therefore gives $P_{\alpha,S}\in\mathcal U(\mathfrak z_L)$ such that
\[
 P_{\alpha,S}|_{E_\beta^\infty}
 =\delta_{\alpha\beta}\operatorname{id}
 \qquad(\beta\in S).
\]

Let $\iota_\alpha:E_\alpha^\infty\to\widehat E$ be the coordinate embedding in \eqref{glob-iso}.
For $w=(w_\beta)\in W$, as $S$ increases the vectors $P_{\alpha,S}w\in W$ converge to $\iota_\alpha(w_\alpha)$ in the product topology.
Closedness gives $\iota_\alpha(w_\alpha)\in W$.
Taking finite partial sums yields
\[
 W=\prod_\alpha W_\alpha,
 \qquad W_\alpha:=\iota_\alpha^{-1}(W).
\]

Each $W_\alpha$ is a closed $L$-subrepresentation of $E_\alpha^\infty$, hence is Casselman--Wallach.
Consequently,
\[
 W_{\mathrm{fin}}=\bigoplus_\alpha(W_\alpha)_{K_L},
 \qquad
 \bigl((W_{\mathrm{fin}})_\alpha\bigr)^\infty\simeq W_\alpha.
\]
As an intersection of $\mathfrak q$-stable subspaces, $W_{\mathrm{fin}}$ is a $\mathcal U(\mathfrak q)$-submodule of $E_{\mathrm{fin}}$.
Lemma~\ref{L:complete-weights} and Noetherianity therefore imply its finite generation.

The same conclusions hold for $W'$.
Thus $W_{\mathrm{fin}}/W'_{\mathrm{fin}}$ is finitely generated and has Harish--Chandra weight spaces, whose globalizations are $W_\alpha/W'_\alpha$ by exactness of globalization.
The product of the coordinate quotient maps is an open surjection with kernel $W'$, giving the asserted topological isomorphism for $W/W'$.
\end{proof}

\section{Bruhat and normal filtrations}\label{S:filtrations}

Retain the principal series $J=C^\infty(X,E)=\operatorname{Ind}_P^G\sigma$ introduced before Corollary~\ref{C:Bruhat-piece-closedness}.
We use the usual tempered structure on this homogeneous bundle.
For a Nash manifold $Y$ and a tempered vector bundle $\mathcal E$, let
\[
 \mathcal S'(Y,\mathcal E):=\mathcal S(Y,\mathcal E)'_b
\]
be continuous functionals on Schwartz test sections.
Here the subscript $b$ denotes the strong dual topology, namely the topology of uniform convergence on bounded subsets of $\mathcal S(Y,\mathcal E)$.
For a closed subset $Z\subset Y$, the subscript $Z$ denotes distributions supported in $Z$.
Since $X$ is compact, $J'=\mathcal S'(X,E)$.

Throughout Sections~\ref{S:filtrations}--\ref{S:comparison-poly}, cohomology means the ordinary cohomology of the underlying Chevalley--Eilenberg complexes.
The topological conclusions will follow by duality in Sections~\ref{S:conclusion} and~\ref{S:Bruhat-piece-closedness}.

Let $\Phi=\Delta(G,A)$ be the restricted root system and $\Phi^+$ be the positive restricted roots defined by $P$.
Put $\Phi_L^+=\Phi(L,A)\cap\Phi^+$, and choose
\[
 {}^LW=\{w\in W(G,A):w^{-1}\Phi_L^+\subset\Phi^+\}.
\]
For each $w\in{}^LW$, choose a representative $\dot w\in N_K(A)$.
The $Q$-orbits are $O_w=Q\dot wP/P$.
Index them compatibly with the ordering in Section~\ref{S:introduction}, so $O_j=O_{w_j}$.
Put $I_k=\mathcal S'_{Z_k}(X,E)$ for $0\leq k\leq r$.
Then $0=I_0\subset\cdots\subset I_r=J'$.
For $O=O_{w_k}$, $1\leq k\leq r$, put
\[
 V_O=X\setminus Z_{k-1},\qquad
 \mathcal T_O=\mathcal S'_O(V_O,E|_{V_O}).
\]
Restriction gives the ordinary Bruhat exact sequence
\begin{equation}\label{Eq:Bruhat-ordinary}
 0\longrightarrow I_{k-1}\longrightarrow I_k
 \longrightarrow\mathcal T_{O_{w_k}}\longrightarrow0.
\end{equation}

\begin{theorem}\label{T:relative-cohomological-comparison}
For $0\leq a\leq b\leq r$, put $D=I_b/I_a$.
For every $i\geq0$, inclusion induces an isomorphism
\[
 H^i\bigl(\mathfrak u,D^{[\mathfrak u]}\bigr)
 \xrightarrow{\ \sim\ }H^i(\mathfrak u,D).
\]
In particular, the assertion holds for $J'=I_r$, for each $\mathcal T_{O_{w_k}}$, and for each $J'/I_j$.
\end{theorem}

Let $N_O=TX|_O/TO$ be the normal bundle.
For $p\geq0$, set
\[
F_p\mathcal T_O
=\left\{T\in\mathcal T_O:
\langle T,\varphi\rangle=0\text{ for all }\varphi\in\mathcal S(V_O,E|_{V_O})
\text{ with }j_O^p\varphi=0\right\},
\qquad F_{-1}\mathcal T_O=0,
\]
where $j_O^p\varphi$ denotes the normal jet of order $p$ along $O$.
Equivalently, $j_O^p\varphi=0$ means that the germ of $\varphi$ belongs to $\mathcal J_O^{p+1}E$, where $\mathcal J_O$ is the ideal of smooth germs vanishing on $O$.
Put
\[
\mathcal E_{p,O}
=E|_O\otimes\operatorname{Sym}^p((N_O\otimes_{\mathbb R}\mathbb C)^\vee),
\qquad
\mathcal G_{p,O}=\mathcal S'(O,\mathcal E_{p,O}).
\]
The normal symbol gives an exact sequence of $Q$-modules
\begin{equation}\label{Eq:normal-ordinary}
 0\longrightarrow F_{p-1}\mathcal T_O
 \longrightarrow F_p\mathcal T_O
 \xrightarrow{\operatorname{symb}_p}\mathcal G_{p,O}
 \longrightarrow0,\qquad
 \mathcal T_O=\bigcup_{p\geq0}F_p\mathcal T_O.
\end{equation}
These facts, including exhaustivity, follow from the Borel lemma; see \cite[Lemma~5.5.3 and Corollary~5.5.4]{AG08}, \cite[Propositions~8.2--8.3]{CS21}, and \cite[Section~2]{Casselman-Hecht-Milicic}.

\begin{lemma}\label{L:comparison-extension}
Let $0\to A\to B\to C\to0$ be an exact sequence of $\mathfrak u$-modules.
If $H^\bullet(\mathfrak u,A^{[\mathfrak u]})\to H^\bullet(\mathfrak u,A)$ is an isomorphism, then
\begin{equation}\label{Eq:extension-exact}
 0\longrightarrow A^{[\mathfrak u]}
 \longrightarrow B^{[\mathfrak u]}
 \longrightarrow C^{[\mathfrak u]}\longrightarrow0
\end{equation}
is exact.
If the comparison also holds for $C$, it holds for $B$.
\end{lemma}

\begin{proof}
Write $\Gamma M=M^{[\mathfrak u]}$.
The comparison for $A$ implies $H^1(\mathfrak u,A/\Gamma A)=0$.
For any $\mathfrak u$-module $M$, we have $H^0(\mathfrak u,M/\Gamma M)=0$: indeed, $\mathfrak u m\subset\Gamma M$ implies $m\in\Gamma M$, since $\mathfrak u$ is finite-dimensional.
Put $D=\operatorname{im}(\Gamma B\to C)\subset\Gamma C$.
The exact sequence
\[
0\longrightarrow A/\Gamma A
\longrightarrow B/\Gamma B
\longrightarrow C/D\longrightarrow0
\]
then gives $H^0(\mathfrak u,C/D)=0$, hence $H^0(\mathfrak u,\Gamma C/D)=0$.
Since $\Gamma C/D$ is locally nilpotent, it must vanish, as every nonzero locally nilpotent $\mathfrak u$-module has a nonzero invariant.
Thus $D=\Gamma C$, proving \eqref{Eq:extension-exact}.
The second assertion follows from the five lemma.
\end{proof}

\begin{proof}[Proof of Theorem~\ref{T:relative-cohomological-comparison} (assuming Proposition~\ref{P:normal-graded-comparison})]
Starting with $F_{-1}\mathcal T_O=0$, apply Lemma~\ref{L:comparison-extension} successively to \eqref{Eq:normal-ordinary}.
Proposition~\ref{P:normal-graded-comparison} gives comparison for every graded piece, so comparison holds for each $F_p\mathcal T_O$.
Since generalized invariants, Chevalley--Eilenberg cochains, and cohomology commute with filtered unions, comparison holds for $\mathcal T_O=\bigcup_pF_p\mathcal T_O$.

We now induct on the Bruhat filtration, starting with $I_0=0$.
If comparison holds for $I_{k-1}$, Lemma~\ref{L:comparison-extension} applied to \eqref{Eq:Bruhat-ordinary} gives an exact sequence
\begin{equation}\label{Eq:Bruhat-exact}
 0\longrightarrow I_{k-1}^{[\mathfrak u]}
 \longrightarrow I_k^{[\mathfrak u]}
 \longrightarrow\mathcal T_{O_{w_k}}^{[\mathfrak u]}
 \longrightarrow0.
\end{equation}
Comparison for the orbit term has already been proved, so the same lemma gives comparison for $I_k$.
Thus comparison holds for every $I_k$, including $I_r=J'$.

Finally, apply Lemma~\ref{L:comparison-extension} to $0\to I_a\to I_b\to D\to0$.
It gives the exact sequence
\[
 0\longrightarrow I_a^{[\mathfrak u]}
 \longrightarrow I_b^{[\mathfrak u]}
 \longrightarrow D^{[\mathfrak u]}\longrightarrow0.
\]
Comparison for $I_a$ and $I_b$, together with the two cohomology long exact sequences, now gives comparison for $D$.
\end{proof}

A direct proof of \eqref{Eq:Bruhat-exact} is given in \cite{Xiao}.

\section{Comparison on a parabolic orbit}\label{S:comparison-poly}

We prove comparison for the normal graded pieces by reducing to scalar distributions on the unipotent fibers.

\begin{proposition}\label{P:normal-graded-comparison}
For every $Q$-orbit $O\subset X$, every $p\geq0$, and every $i\geq0$, inclusion induces an isomorphism
\begin{equation}\label{Eq:parameter-comparison}
 H^i(\mathfrak u,\mathcal G_{p,O}^{[\mathfrak u]})
 \xrightarrow{\ \sim\ }H^i(\mathfrak u,\mathcal G_{p,O}).
\end{equation}
\end{proposition}

\subsection{Reduction to the scalar case}

Fix $w\in{}^LW$, put $x_w=\dot wP\in X$, and write $O=O_w=Qx_w$.
Let $\bar N$ be the unipotent radical of the parabolic opposite to $P$ with common Levi subgroup $MA$.
Put
\[
P_L=L\cap P,\qquad H_w=U\cap\dot wN\dot w^{-1},\qquad
U_w^-=U\cap\dot w\bar N\dot w^{-1},\qquad Y_w=U/H_w.
\]

\begin{lemma}\label{L:orbit-fibration}
The stabilizer of $x_w$ in $Q$ is $S_w=P_L\ltimes H_w$, and multiplication induces an algebraic isomorphism $U_w^-\times H_w\simeq U$.
The group $P_L$ acts on $Y_w$ by conjugation, and there is an algebraic isomorphism
\begin{equation}\label{Eq:orbit-fibration}
O\simeq Q/S_w\simeq L\times^{P_L}Y_w
\xrightarrow{\ \pi_w\ }B:=L/P_L,
\qquad [l,uH_w]\longmapsto lP_L,
\end{equation}
where $(la,y)\sim(l,aya^{-1})$ for $l\in L$, $a\in P_L$, and $y\in Y_w$, and $[l,uH_w]$ corresponds to $lu x_w$ for $u\in U$.
The base $B$ is compact and $U$ acts trivially on it.
\end{lemma}

\begin{proof}
Write $P^w=\dot wP\dot w^{-1}$, so that $S_w=Q\cap P^w$.
The condition $w^{-1}\Phi_L^+\subset\Phi^+$ gives $L\cap P^w=P_L$.
Also $U\cap P^w=H_w$, since the Levi factor $MA$ of $P^w$ has no nontrivial unipotent elements.

Choose $H_Q\in\operatorname{Lie}(A)$ centralizing $L$ with positive integral weights on $\operatorname{Lie}(U)$.
Then $a_t=\exp(-tH_Q)$ centralizes $L$ and contracts $U$ as $t\to\infty$.
If $lu\in S_w$, then $a_t(lu)a_t^{-1}\to l$.
Since $a_t\in P^w$ and $P^w$ is closed, we get $l\in P_L$ and hence $u\in H_w$.
Thus $S_w=P_L\ltimes H_w$.

In particular, $P_L$ normalizes $H_w$ and acts on $Y_w=U/H_w$ by conjugation.
Then
\[
 L\times^{P_L}Y_w \xrightarrow{\ \sim\ } Q/S_w=O,
 \qquad [l,uH_w]\longmapsto luS_w.
\]
Under this identification, $\pi_w$ is induced by the Levi projection $Q\to L$, so $U$ acts trivially on the base.
The base $L/P_L$ is compact because $P_L$ is a minimal parabolic subgroup of $L$.

The multiplication assertion follows from \cite[Proposition~3.11]{BT65} by partitioning the restricted roots $\alpha$ occurring in $\operatorname{Lie}(U)$ according to the sign of $w^{-1}\alpha$.
\end{proof}

Define the algebraic section $c_w:Y_w\to U_w^-\subset U$ by $c_w(vH_w)=v$ for $v\in U_w^-$.
Choose $l_i\in L$ such that the translated opposite big cells $B_i:=l_i(L\cap\bar N)P_L/P_L$ cover $B=L/P_L$, and define sections $s_i:B_i\to L$ by $s_i(l_i vP_L)=l_i v$ for $v\in L\cap\bar N$.

The following algebraic chart and its lift are fixed throughout:
\begin{equation}\label{Eq:orbit-chart}
\begin{aligned}
\theta_i:B_i\times Y_w&\xrightarrow{\ \sim\ }\Omega_i:=\pi_w^{-1}(B_i),
&\qquad \theta_i(b,y)&=s_i(b)c_w(y)x_w,\\
q_i:B_i\times Y_w&\longrightarrow Q,
&\qquad q_i(b,y)&=s_i(b)c_w(y).
\end{aligned}
\end{equation}
Thus $q_i$ lifts $\theta_i$ through the orbit map $Q\to O$, $q\mapsto qx_w$.
The $U$-action on $B_i\times Y_w$ is
\begin{equation}\label{Eq:twisted-action}
 n\cdot(b,y)=(b,a_b(n)y),\qquad a_b(n)=s_i(b)^{-1}ns_i(b).
\end{equation}

For a compact $K\subset B_i$, set
\[
 M_{p,i,K}=\mathcal S'_{\pi_w^{-1}(K)}(\Omega_i,\mathcal E_{p,O}),
 \qquad \mathcal T_K(Y_w)=\mathcal S'_{K\times Y_w}(B_i\times Y_w).
\]
We will consider two $U$-actions on $\mathcal T_K(Y_w)$, induced respectively by the fixed fiber action $n\cdot(b,y)=(b,ny)$ and the natural twisted action \eqref{Eq:twisted-action}.

The following lemma is immediate and independent of the cutoff.

\begin{lemma}\label{L:orbit-support-extension}
Choose $\eta_i\in C_c^\infty(B_i)$ equal to one near $K$.
The formula
\begin{equation}\label{Eq:orbit-support-extension}
 \langle e_iT,\phi\rangle
 =\bigl\langle T,
       ((\pi_w^*\eta_i)\phi)|_{\Omega_i}\bigr\rangle
 \qquad(\phi\in\mathcal S(O,\mathcal E_{p,O}))
\end{equation}
defines a continuous $U$-equivariant extension $e_i:M_{p,i,K}\to\mathcal G_{p,O}$, supported in $\pi_w^{-1}(K)$.
It preserves $[\mathfrak u]$-invariants.
\end{lemma}

\begin{lemma}\label{L:normal-coefficient-filtration}
The module $M_{p,i,K}$ admits a finite filtration by $\mathfrak u$-submodules whose successive quotients are $\mathcal T_K(Y_w)\otimes W_j^*$ for some finite-dimensional vector spaces $W_j$, where $\mathfrak u$ acts on $\mathcal T_K(Y_w)$ through \eqref{Eq:twisted-action} and trivially on $W_j^*$.
\end{lemma}

\begin{proof}
The normal space at $x_w$ and its dual symmetric power
\[
 N_{O,x_w,\mathbb C}\simeq
 \mathfrak g/(\mathfrak q+\operatorname{Ad}(\dot w)\mathfrak p),
 \qquad V_p=\operatorname{Sym}^pN_{O,x_w,\mathbb C}^*
\]
are algebraic $S_w$-modules.
$S_w$ acts naturally on the fiber $V_\sigma^w:=E_{x_w}$, and by identifying $V_\sigma^w$ with $V_\sigma$ via $v\mapsto[\dot w,v]$, this action is given by $s\mapsto\sigma(\dot w^{-1}s\dot w)$ for $s\in S_w$.
In particular, $H_w$ acts trivially on $V_\sigma^w$.
We have a canonical isomorphism of $Q$-equivariant tempered bundles
\[
 \mathcal E_{p,O}\simeq Q\times^{S_w}(V_\sigma^w\otimes V_p).
\]
The spaces
\[
 V_{p,j}=\ker(I_{\mathfrak h_w}^{j};V_p),\qquad
 0=V_{p,0}\subset\cdots\subset V_{p,m}=V_p
\]
form a finite $S_w$-stable filtration of $V_p$.
Each quotient is $H_w$-trivial.
Hence the subbundles
\[
 \mathcal E_j=Q\times^{S_w}(V_\sigma^w\otimes V_{p,j})
 \subset\mathcal E_{p,O}
\]
have quotients $\mathcal F_j:=\mathcal E_j/\mathcal E_{j-1}\simeq Q\times^{S_w}W_j$, where $W_j=V_\sigma^w\otimes(V_{p,j}/V_{p,j-1})$.
Since $H_w$ acts trivially on $W_j$, the projection $r:Q=L\ltimes U\to L$ induces the $Q$-equivariant isomorphism
\[
\begin{aligned}
 Q\times^{S_w}W_j&\xrightarrow{\ \sim\ }\pi_w^*(L\times^{P_L}W_j),\\
 [q,v]&\longmapsto\bigl(qS_w,[r(q),v]\bigr).
\end{aligned}
\]

Here the $Q$-action on the pullback bundle is given by $q\cdot(x,\xi)=(qx,r(q)\cdot\xi)$.
Choose a basis $(e_{j,a})_a$ of $W_j$.
The $U$-invariant sections
\[
 \tau_{i,j,a}(\theta_i(b,y))=[q_i(b,y),e_{j,a}]
 \longleftrightarrow\bigl(\theta_i(b,y),[s_i(b),e_{j,a}]\bigr)
\]
give a local frame of $\mathcal F_j|_{\Omega_i}$, with trivialization
\[
 (B_i\times Y_w)\times W_j\xrightarrow{\ \sim\ }\mathcal F_j|_{\Omega_i},
 \qquad (b,y,v)\longmapsto[q_i(b,y),v].
\]
The action in this trivialization is $n\cdot(b,y,v)=\bigl(b,a_b(n)y,v\bigr)$.
In this frame, $\mathcal S(\Omega_i,\mathcal F_j)\simeq\mathcal S(B_i\times Y_w)\otimes W_j$.
Passing to continuous duals with the prescribed support gives
\[
 \mathcal S'_{\pi_w^{-1}(K)}(\Omega_i,\mathcal F_j)
 \simeq\mathcal T_K(Y_w)\otimes W_j^*,
\]
where $\mathfrak u$ acts through \eqref{Eq:twisted-action} on the first factor and trivially on the second.

Define the annihilator filtration by
\[
 M_j=\{T\in M_{p,i,K}:
       T|_{\mathcal S(\Omega_i,\mathcal E_{m-j})}=0\},
 \qquad 0\leq j\leq m.
\]
We obtain a filtration by $\mathfrak u$-submodules with
\[
 0=M_0\subset\cdots\subset M_m=M_{p,i,K},\qquad
 M_j/M_{j-1}\simeq\mathcal T_K(Y_w)\otimes W_{m-j+1}^*.
\]
\end{proof}

We use the following scalar comparison, whose proof is deferred to Subsection~\ref{S:scalar-comparison}.

\begin{lemma}\label{L:scalar-parameter-comparison}
Let $N_1$ be a real unipotent algebraic group, $H_1$ an algebraic subgroup, $Y=N_1/H_1$, and $\mathfrak n_1=\operatorname{Lie}(N_1)_{\mathbb C}$.
For a Nash manifold $D$ and a compact subset $K\subset D$, put $\mathcal T_K(Y)=\mathcal S'_{K\times Y}(D\times Y)$, with $N_1$ acting trivially on $D$ and by left translation on $Y$.
Then inclusion induces isomorphisms
\[
 H^i(\mathfrak n_1,\mathcal T_K(Y)^{[\mathfrak n_1]})
 \xrightarrow{\ \sim\ }H^i(\mathfrak n_1,\mathcal T_K(Y))
 \qquad(i\geq0).
\]
\end{lemma}

\begin{proof}[Proof of Proposition~\ref{P:normal-graded-comparison}]
Put $M=\mathcal T_K(Y_w)$.
We first transfer the scalar comparison from the translation action to the twisted action \eqref{Eq:twisted-action}.
Write $C^\bullet_{\mathrm{tr}}(M)$ and $C^\bullet_{\mathrm{tw}}(M)$ for the corresponding Chevalley--Eilenberg complexes, and put $a_b=\operatorname{Ad}(s_i(b)^{-1})|_{\mathfrak u}$.
Define $\Phi:C^\bullet_{\mathrm{tw}}(M)\to C^\bullet_{\mathrm{tr}}(M)$ by
\begin{equation}\label{Eq:cochain-untwisting}
 (\Phi c)(X_1,\ldots,X_q)
       =c(a_b^{-1}X_1,\ldots,a_b^{-1}X_q).
\end{equation}
The matrix coefficients of $a_b^{\pm1}$ depend only on $b$, so multiplication by them commutes with the fiber derivatives.
Hence $\Phi$ is a cochain isomorphism.

For the same reason, every product of $r$ twisted operators is a linear combination, with coefficients depending only on $b$, of products of $r$ translation operators, and conversely.
Thus the two actions have the same $I_{\mathfrak u}^r$-annihilated subspace, which is preserved by $\Phi$ coefficientwise.
Thus we have a commutative diagram
\[
\begin{tikzcd}[column sep=large, row sep=large]
 C^\bullet_{\mathrm{tw}}(M^{[\mathfrak u]}) \arrow[r,hook] \arrow[d,"\Phi"']
 & C^\bullet_{\mathrm{tw}}(M) \arrow[d,"\Phi"] \\
 C^\bullet_{\mathrm{tr}}(M^{[\mathfrak u]}) \arrow[r,hook]
 & C^\bullet_{\mathrm{tr}}(M),
\end{tikzcd}
\]
whose vertical arrows are isomorphisms of complexes.
By Lemma~\ref{L:scalar-parameter-comparison} with $N_1=U$ and $H_1=H_w$, the lower horizontal arrow is a quasi-isomorphism, and hence so is the upper one.
Tensoring with each finite-dimensional $W_j^*$ gives comparison for every quotient in Lemma~\ref{L:normal-coefficient-filtration}.
Induction along its finite filtration, using Lemma~\ref{L:comparison-extension}, now yields
\begin{equation}\label{Eq:local-normal-comparison}
 H^q(\mathfrak u,M_{p,i,K}^{[\mathfrak u]})
 \xrightarrow{\ \sim\ }H^q(\mathfrak u,M_{p,i,K}).
\end{equation}

Choose a finite smooth partition of unity $\sum_i\beta_i=1$ on $B$ with $K_i=\operatorname{supp}\beta_i\Subset B_i$.
Multiplication by $\pi_w^*\beta_i$ commutes with $\mathfrak u$ and with the cochain differential $d$.
For a cocycle $z\in C^q(\mathfrak u,\mathcal G_{p,O})$, local surjectivity in \eqref{Eq:local-normal-comparison} gives
\[
 (\beta_i z)|_{\Omega_i}=p_i+db_i,
 \quad p_i\in Z^q(\mathfrak u,M_{p,i,K_i}^{[\mathfrak u]}),
 \quad b_i\in C^{q-1}(\mathfrak u,M_{p,i,K_i}).
\]
Lemma~\ref{L:orbit-support-extension} therefore yields $z=\sum_i e_ip_i+d\sum_i e_ib_i$, proving surjectivity.
If a cocycle $p$ with generalized-invariant coefficients satisfies $p=db$ in the full complex, then $(\beta_ip)|_{\Omega_i}=d((\beta_ib)|_{\Omega_i})$.
Local injectivity supplies $a_i\in C^{q-1}(\mathfrak u,M_{p,i,K_i}^{[\mathfrak u]})$ with $da_i=(\beta_ip)|_{\Omega_i}$.
Hence $p=d\sum_i e_ia_i$.
This proves \eqref{Eq:parameter-comparison}.
\end{proof}

\subsection{The scalar comparison}\label{S:scalar-comparison}

\begin{lemma}\label{L:polynomial-parameters}
In the setting of Lemma~\ref{L:scalar-parameter-comparison}, let $\delta_Y$ be the density associated with a nonzero $N_1$-invariant algebraic top form on $Y$, and let $\mathcal D'_K(D)$ denote the distributions on $D$ supported in $K$.
Then
\begin{equation}\label{Eq:polynomial-parameters}
 \mathcal T_K(Y)^{[\mathfrak n_1]}
 =\mathcal D'_K(D)\otimes R(Y)\delta_Y
 =:\mathcal P_K(Y),
\end{equation}
where the tensor product is algebraic and is embedded in $\mathcal T_K(Y)$ by
\begin{equation}
 \langle T\otimes f\delta_Y,\phi\rangle
 =\left\langle T,\ b\longmapsto
          \int_Y f(y)\phi(b,y)\,\delta_Y(y)\right\rangle.
\end{equation}
\end{lemma}

\begin{proof}
Write $I=I_{\mathfrak n_1}$ and fix $r\geq1$.
By \cite[Lemma~3.3(i)]{LLY}, the scalar distributions on $Y$ annihilated by $I^r$ form the finite-dimensional space $E_r=(R(Y)\delta_Y)^{I^r}$.
If $T\in\mathcal T_K(Y)$ satisfies $I^rT=0$, then for every $\phi\in\mathcal S(D)$ its slice
\[
 T_\phi(\psi):=\langle T,\phi\otimes\psi\rangle
 \qquad(\psi\in\mathcal S(Y))
\]
belongs to $E_r$, since $N_1$ acts only on $Y$.

Choose a basis $e_1,\ldots,e_m$ of $E_r$ and test functions $\psi_a\in\mathcal S(Y)$ with $\langle e_j,\psi_a\rangle=\delta_{ja}$.
Expanding each slice in this basis gives
\[
 T_\phi=\sum_{a=1}^m T_a(\phi)e_a,\qquad
 T_a(\phi)=\langle T,\phi\otimes\psi_a\rangle .
\]
Each $T_a$ is continuous on $\mathcal S(D)$ and vanishes when $\phi$ vanishes near $K$, so $T_a\in\mathcal D'_K(D)$.
By density of tensor test functions in $\mathcal S(D\times Y)$, we obtain $T=\sum_a T_a\otimes e_a$.

Conversely, $I^r$ annihilates $\mathcal D'_K(D)\otimes E_r$ because it acts only on the second factor.
Thus $\mathcal T_K(Y)^{I^r}=\mathcal D'_K(D)\otimes E_r$.
By \cite[Lemma~3.3(ii)]{LLY}, $\bigcup_r E_r=R(Y)\delta_Y$.
Taking the union over $r$ gives \eqref{Eq:polynomial-parameters}, since every algebraic tensor is a finite sum.
\end{proof}

\begin{proof}[Proof of Lemma~\ref{L:scalar-parameter-comparison}]
By Lemma~\ref{L:polynomial-parameters}, it suffices to show that the inclusion $\mathcal P_K(Y)\hookrightarrow\mathcal T_K(Y)$ induces an isomorphism on $\mathfrak n_1$-cohomology.
Following \cite[Section~4]{LLY}, we induct on $\dim N_1$ with $D$ and $K$ fixed.
The case $N_1=1$ is immediate.
Choose a one-dimensional central algebraic subgroup $C\subset N_1$ and put $\mathfrak c=\operatorname{Lie}(C)_{\mathbb C}$.
For $M=\mathcal P_K(Y)$ or $\mathcal T_K(Y)$, consider the Hochschild--Serre spectral sequence
\[
 E_2^{a,b}(M)=H^a(\mathfrak n_1/\mathfrak c,
                        H^b(\mathfrak c,M))
 \ \Longrightarrow\ H^{a+b}(\mathfrak n_1,M).
\]

\emph{If $C\subset H_1$,} then $C$ acts trivially on $Y$, so $H^b(\mathfrak c,M)=\bigwedge^b\mathfrak c^*\otimes M$.
The induction hypothesis for $(N_1/C,H_1/C)$ identifies the two $E_2$ pages.
Spectral-sequence comparison proves the assertion.

\emph{If $C\not\subset H_1$,} then $C\cap H_1=1$, and the fibers of $p:Y\to\overline Y=N_1/(CH_1)$ are free $C$-orbits.
To compute the $\mathfrak c$-cohomology, we choose a coordinate along the $C$-orbits.
Put $B=CH_1$ and fix $0\ne Z\in\operatorname{Lie}(C)$.
By \cite[Theorem~1.1.13(a) and Proposition~1.2.8]{CG90}, ordered exponential coordinates adapted to $\operatorname{Lie}(B)$ give an algebraic section $\sigma:\overline Y=N_1/B\to N_1$ of the quotient map, obtained by setting the $B$-coordinates equal to zero.

Since $B=C\times H_1$, this gives an algebraic isomorphism
\[
 \mathbb A^1_{\mathbb R}\times\overline Y\xrightarrow{\ \sim\ }Y,
 \qquad (t,\bar y)\longmapsto\exp(tZ)\sigma(\bar y)H_1.
\]
In these coordinates $p$ is projection onto $\overline Y$, and $\exp(aZ)$ acts by $t\mapsto t+a$.
Thus $R(Y)=R(\overline Y)[t]$, and $Z$ acts on test functions by $-\partial_t$.

Normalize the invariant densities so that $\delta_Y=|dt|\,\delta_{\overline Y}$; see also \cite[Lemma~4.1]{LLY}.
On $D\times Y$ the same coordinates leave the parameter in $D$ unchanged, so integration and differentiation in $t$ preserve the prescribed base support.

Choose $\rho\in\mathcal S(\mathbb R)$ with $\int\rho=1$.
Define continuous operators
\[
\begin{aligned}
 &A:\mathcal S(D\times Y)\longrightarrow\mathcal S(D\times\overline Y),\qquad
   A\phi(b,\bar y)=\int_{\mathbb R}\phi(b,t,\bar y)\,dt,\\
 &j:\mathcal S(D\times\overline Y)\longrightarrow\mathcal S(D\times Y),\qquad
   j\psi(b,t,\bar y)=\rho(t)\psi(b,\bar y),\\
 &h:\mathcal S(D\times Y)\longrightarrow\mathcal S(D\times Y),
 \qquad h\phi(b,t,\bar y)=\int_{-\infty}^{t}\bigl(\phi(b,s,\bar y)-\rho(s)A\phi(b,\bar y)\bigr)\,ds.
\end{aligned}
\]

Since the operators act only in the $t$ variable, their transposes preserve base support in $K$ and restrict to
\[
\begin{aligned}
 A':\mathcal T_K(\overline Y)\longrightarrow\mathcal T_K(Y),
 j':\mathcal T_K(Y)\longrightarrow\mathcal T_K(\overline Y),
 h':\mathcal T_K(Y)\longrightarrow\mathcal T_K(Y).
\end{aligned}
\]
The identities $Aj=1,\ A\partial_t=0,\ \partial_t h=1-jA,\ h\partial_t=1$ give
\[
 j'A'=1,\qquad \partial_t'A'=0,\qquad \partial_t'h'=1,\qquad h'\partial_t'=1-A'j'.
\]
Thus $\partial_t':\mathcal T_K(Y)\to\mathcal T_K(Y)$ is surjective, $A'$ is injective, and $\ker\bigl(\partial_t':\mathcal T_K(Y)\to\mathcal T_K(Y)\bigr)=A'\mathcal T_K(\overline Y)$.
For $\mathcal P_K$, use $R(Y)=R(\overline Y)[t]$.
It follows that $\partial_t':\mathcal P_K(Y)\to\mathcal P_K(Y)$ is also surjective, with kernel $A'\mathcal P_K(\overline Y)$.
Since $\mathfrak c$ is one-dimensional and its Chevalley--Eilenberg differential is $\partial_t'$, we obtain
\[
 H^b(\mathfrak c,\mathcal T_K(Y))=
 \begin{cases}\mathcal T_K(\overline Y),&b=0,\\0,&b>0,\end{cases}
\]
and likewise for $\mathcal P_K$, compatibly with inclusion.
The degree-zero identification by $A'$ is $N_1/C$-equivariant, hence induction for $(N_1/C,CH_1/C)$ identifies the $E_2$ pages, and the spectral sequences prove the comparison.
This completes the proof.
\end{proof}

\section{Proof of the comparison theorem}\label{S:conclusion}

For $\pi\in\mathrm{CW}_G$, recall that
\[
 \mathrm{CJ}^{\infty}_U(\pi)=\varprojlim_{k\geq1}\pi/\overline{\mathfrak u^k\pi},
 \qquad
 q_\pi: \pi\to \mathrm{CJ}^{\infty}_U(\pi),\  v\mapsto\bigl(v+\overline{\mathfrak u^k\pi}\bigr)_{k\geq1}.
\]
Every continuous linear functional on this inverse limit factors through a finite-level quotient, giving the identification
\begin{equation}\label{dual iso}
 (\mathrm{CJ}^{\infty}_U(\pi))'
 =\bigcup_{k\geq1}(\pi/\overline{\mathfrak u^k\pi})'
 =(\pi')^{[\mathfrak u]}.
\end{equation}
Under this identification, the transpose $q_\pi'$ is the inclusion $(\pi')^{[\mathfrak u]}\hookrightarrow\pi'$.

\begin{theorem}\label{T:smooth-CJ-comparison}
For every Casselman--Wallach representation $\pi$ of $G$, the quotient map $q_\pi$ induces topological isomorphisms
\begin{equation}\label{CW-comp}
 H_i(\mathfrak u,\pi)\xrightarrow{\ \sim\ }
 H_i(\mathfrak u,\mathrm{CJ}^{\infty}_U(\pi)),\qquad i\geq0.
\end{equation}
\end{theorem}

\begin{lemma}\label{L:Frechet-duality}
Let $f:C_\bullet\to D_\bullet$ be a continuous map of bounded complexes of Fr\'echet spaces.
\begin{enumerate}
\item If the transpose $f':D_\bullet'\to C_\bullet'$ is a quasi-isomorphism, then $f$ is a quasi-isomorphism of the underlying complexes of vector spaces.
\item If $f$ is a quasi-isomorphism of the underlying complexes and the homology of $D_\bullet$ is Hausdorff, then the homology of $C_\bullet$ is Hausdorff and the induced isomorphisms are topological.
\end{enumerate}
\end{lemma}

\begin{proof}
For (1), let $(K_\bullet,d)=\operatorname{Cone}(f)$.
Its dual is exact, so $\operatorname{im}d_i'=\ker d_{i+1}'$ is weak-* closed.
By the closed range theorem for Fr\'echet spaces, every $d_i$ has closed image.
Taking preannihilators and using Hahn--Banach gives
\[
 \operatorname{im}d_{i+1}={}^{\perp}\ker d_{i+1}'={}^{\perp}\operatorname{im}d_i'=\ker d_i.
\]
Thus $K_\bullet$ is exact and $f$ is a quasi-isomorphism.

For (2), put $Z_i(C)=\ker d_i^C$ and $B_i(C)=\operatorname{im}d_{i+1}^C$.
The continuous map $Z_i(C)\to H_i(D_\bullet)$ induced by $f$ has kernel $B_i(C)$, because $H_i(f)$ is injective.
Hence $B_i(C)$ is closed in $Z_i(C)$, and therefore in $C_i$, since $Z_i(C)$ is closed.
Both $H_i(C_\bullet)$ and $H_i(D_\bullet)$ are consequently Fr\'echet spaces.
The map $H_i(f)$ is a continuous bijection by hypothesis, so the open mapping theorem makes it a topological isomorphism.
\end{proof}

Every finite-dimensional smooth $P$-module $F$ admits a finite $P$-stable filtration whose successive quotients have the form $\sigma_M\otimes e^\lambda$, with $N$ acting trivially.
Smooth induction gives strict exact sequences of Casselman--Wallach representations, which remain strictly exact after applying $\mathrm{CJ}^{\infty}_U$ by \cite[Theorem~14.3]{CWYZ}.

Taking continuous duals and using \eqref{dual iso}, induction along the filtration extends the comparison for the full dual $J'$ in Theorem~\ref{T:relative-cohomological-comparison} to $J=\operatorname{Ind}_P^G F$.
Applying parts (1) and (2) of Lemma~\ref{L:Frechet-duality}, with the target homology Hausdorff by \eqref{Eq:CJ-comparison}, gives the following corollary.

\begin{corollary}\label{P:principal-CJ-comparison}
Let $F$ be a finite-dimensional smooth representation of the fixed minimal parabolic subgroup $P$, and put $J=\operatorname{Ind}_P^G F$.
Then $q_J$ induces topological isomorphisms
\[
 H_i(\mathfrak u,J)\xrightarrow{\ \sim\ }
 H_i(\mathfrak u,\mathrm{CJ}^{\infty}_U(J)),\qquad i\geq0.
\]
\end{corollary}

\begin{proof}[Proof of Theorem~\ref{T:smooth-CJ-comparison}]
Casselman's embedding theorem gives a closed embedding of every Casselman--Wallach representation into $\operatorname{Ind}_P^G F$ for a finite-dimensional smooth $P$-representation $F$.
Applying it successively to the Casselman--Wallach cokernels gives a strict coresolution
\[
 0\longrightarrow\pi\longrightarrow J^0\longrightarrow J^1
 \longrightarrow\cdots,\qquad J^p=\operatorname{Ind}_P^G F^p;
\]
see \cite[Proposition~4.2.3]{Wallach1} and \cite[Proposition~3]{Hecht-Taylor3}.
By \cite[Theorem~14.3]{CWYZ}, strict exactness of $\mathrm{CJ}^{\infty}_U$ gives a strict coresolution of $\mathrm{CJ}^{\infty}_U(\pi)$ by $\mathrm{CJ}^{\infty}_U(J^p)$.

Form the two double complexes with terms
\[
 J^p\otimes\bigwedge^q\mathfrak u,
 \qquad
 \mathrm{CJ}^{\infty}_U(J^p)\otimes\bigwedge^q\mathfrak u,
 \qquad p\geq0,\quad 0\leq q\leq\dim\mathfrak u.
\]
Use total homological degree $q-p$.
Since $q$ ranges over a finite interval, each total degree has only finitely many terms, so the two spectral sequences converge.
The maps $q_{J^p}$ give a morphism of double complexes.
The coresolutions identify the total homology with the homology of $\pi$ and $\mathrm{CJ}^{\infty}_U(\pi)$, respectively.
By Corollary~\ref{P:principal-CJ-comparison}, the induced maps on column homology are isomorphisms.
Spectral-sequence comparison therefore gives
\begin{equation}\label{CW-com-prove}
 H_i(\mathfrak u,\pi)\simeq
 H_i(\mathfrak u,\mathrm{CJ}^{\infty}_U(\pi)).
\end{equation}
This is the map induced by the continuous chain map $C_\bullet(\mathfrak u,q_\pi)$.
Its target is Casselman--Wallach, and hence Hausdorff, by \eqref{Eq:CJ-comparison}.
Part (2) of Lemma~\ref{L:Frechet-duality} therefore shows that \eqref{CW-com-prove} is a topological isomorphism and that the source homology is Hausdorff.
In particular, the differentials of $C_\bullet(\mathfrak u,\pi)$ have closed images.
\end{proof}

\begin{proof}[Proof of Theorem~\ref{comparison conj}]
Corollary~\ref{C:CJ-comparison} and Theorem~\ref{T:smooth-CJ-comparison} give natural topological $L$-equivariant isomorphisms
\begin{equation}\label{iso comp eq}
 H_i(\mathfrak u,\pi_K)^\infty
 \xrightarrow{\ \sim\ }
 H_i(\mathfrak u,\mathrm{CJ}^{\infty}_U(\pi))
 \xleftarrow{\ \sim\ }H_i(\mathfrak u,\pi),
\end{equation}
and the composition extends the map induced by $\pi_K\hookrightarrow\pi$.
Theorem~\ref{T:smooth-CJ-comparison} also gives closed images for the differentials of $C_\bullet(\mathfrak u,\pi)$.
\end{proof}

\begin{proof}[Proof of Corollary~\ref{C:cohomological-comparison}]
Put $d=\dim_{\mathbb C}\mathfrak u$ and $D=(\bigwedge^d\mathfrak u)^*$.
Since $\mathfrak u$ is nilpotent, $\mathfrak u$ acts trivially on $D$.
For $W=\pi_K$ or $\pi$, we have a natural isomorphism of complexes
\[
 C^\bullet(\mathfrak u,W)
 \simeq C_{d-\bullet}(\mathfrak u,W)\otimes D.
\]
It is $(\mathfrak l,K_L)$-equivariant for $W=\pi_K$ and $L$-equivariant for $W=\pi$.
In the latter case it is a homeomorphism in each degree, hence also for the induced quotient topologies on cohomology.

The cohomological comparison map in degree $i$ identifies with the homological comparison map in degree $d-i$ tensored with $\mathrm{id}_D$, hence is a topological $L$-isomorphism by Theorem~\ref{comparison conj}.
The same isomorphism of complexes transfers closed images to the cochain complex.
\end{proof}

\section{Finite-dimensional coefficients and closedness}\label{S:finite-coefficients}

We first extend comparison to finite-dimensional coefficients and deduce Corollaries~\ref{C:augmentation-comparison} and~\ref{C:Koszul-component-closedness}.
For Corollary~\ref{C:Bruhat-piece-closedness}, we first establish closed differential images and Casselman--Wallach homology on the smooth Bruhat pieces.

\subsection{Finite-dimensional coefficients}\label{S:finite-coefficient-comparison}

We use the following elementary closed-image argument.
\begin{lemma}\label{L:closed-range-extension}
Let $0\to A_\bullet\to B_\bullet\to C_\bullet\to0$ be a degreewise strict exact sequence of bounded complexes of smooth moderate-growth Fr\'echet $L$-representations.
If the differentials of $A_\bullet$ and $C_\bullet$ have closed images and their homology belongs to $\mathrm{CW}_L$, the same holds for $B_\bullet$.
\end{lemma}
\begin{proof}
The connecting maps $\delta_i:H_i(C_\bullet)\to H_{i-1}(A_\bullet)$ are continuous by the open mapping theorem.
They are therefore morphisms in $\mathrm{CW}_L$ and have closed images.

Write $Z_i(A)=\ker d_i^A$, and similarly for $B$.
The space
\[
 K_i=\ker\bigl(Z_i(B)\to H_i(C_\bullet)\bigr)
     =Z_i(A)+dB_{i+1}
\]
is closed in $B_i$.
The continuous surjection
\[
 Z_i(A)\oplus B_{i+1}\longrightarrow K_i,
 \qquad(a,b)\longmapsto a+db
\]
is open.
The inverse image of $dB_{i+1}$ is $W_i\oplus B_{i+1}$, where $W_i$ is the inverse image of $\operatorname{im}\delta_{i+1}$ under $Z_i(A)\to H_i(A_\bullet)$.
Thus $W_i$ is closed, and hence $dB_{i+1}$ is closed in $K_i$ and in $B_i$.

The homology long exact sequence and the open mapping theorem now give a strict exact sequence
\[
 0\longrightarrow\operatorname{coker}\delta_{i+1}
 \longrightarrow H_i(B_\bullet)
 \longrightarrow\ker\delta_i\longrightarrow0.
\]
Its middle term is a smooth moderate-growth Fr\'echet representation, and its $K_L$-finite vectors form an extension of Harish--Chandra modules.
Consequently $H_i(B_\bullet)$ belongs to $\mathrm{CW}_L$.
\end{proof}

\begin{proposition}\label{P:finite-coefficient-comparison}
Let $F$ be a finite-dimensional algebraic $Q$-module and let $\pi\in\mathrm{CW}_G$.
For every $i\geq0$, $H_i(\mathfrak u,F\otimes\pi_K)$ is Harish--Chandra, and inclusion induces a topological $L$-isomorphism
\begin{equation}\label{Eq:finite-coefficient-comparison}
 H_i(\mathfrak u,F\otimes\pi_K)^\infty
 \xrightarrow{\ \sim\ }H_i(\mathfrak u,F\otimes\pi),
\end{equation}
natural in $F$ and $\pi$.
The smooth complex has closed differential images and Casselman--Wallach homology.
The corresponding assertions also hold for cohomology.

More generally, if $M$ is a smooth moderate-growth Fr\'echet $(\mathfrak u,L)$-module whose Chevalley--Eilenberg complex has closed differential images and Casselman--Wallach homology, these properties also hold for $F\otimes M$.
\end{proposition}
\begin{proof}
Since $U$ is unipotent, $F$ has a finite $L$-stable filtration
\[
 F\supset I_{\mathfrak u}F\supset\cdots\supset I_{\mathfrak u}^sF=0
\]
whose successive quotients have trivial $\mathfrak u$-action.
For such a quotient $A$, one has
\[
 C_\bullet(\mathfrak u,A\otimes M)
 =A\otimes C_\bullet(\mathfrak u,M).
\]
Lemma~\ref{L:closed-range-extension}, applied along this filtration, shows that $C_\bullet(\mathfrak u,F\otimes M)$ has closed differential images and Casselman--Wallach homology.
Taking $M=\pi$, for which Theorem~\ref{comparison conj} supplies the hypotheses, gives the same conclusions for $C_\bullet(\mathfrak u,F\otimes\pi)$.

On the algebraic side, for each successive quotient $A$ one has
\[
 H_i(\mathfrak u,A\otimes\pi_K)
 =A\otimes H_i(\mathfrak u,\pi_K),
\]
which is Harish--Chandra by Theorem~\ref{comparison conj}, so $H_i(\mathfrak u,F\otimes\pi_K)$ is also Harish--Chandra.

The inclusion induces an $(\mathfrak l,K_L)$-equivariant map
\[
 H_i(\mathfrak u,F\otimes\pi_K)
 \longrightarrow H_i(\mathfrak u,F\otimes\pi).
\]
Its source is now Harish--Chandra and its target Casselman--Wallach, so it extends uniquely to a continuous $L$-map from the globalization of the source.

For each successive quotient $A$, this extension is an isomorphism by Theorem~\ref{comparison conj} tensored with $A$.
Exactness of globalization allows us to compare the two homology long exact sequences along the filtration; the five lemma then gives \eqref{Eq:finite-coefficient-comparison}.
Naturality of inclusion and uniqueness of extension give naturality in $F$ and $\pi$.

The cohomological assertions follow from the duality argument in the proof of Corollary~\ref{C:cohomological-comparison}.
\end{proof}

Equip $F_k$ with (left) adjoint action of $L$.
Regard $F_k$ first as a right $\mathscr U$-module, and write $F_k^-$ for the left $\mathfrak u$-module with action $X\cdot a=-aX$.
These actions make $F_k^-$ a finite-dimensional algebraic (left) $Q$-module.

For a left $\mathscr U$-module $V$, the standard Chevalley--Eilenberg resolution computes $\operatorname{Tor}_i^{\mathscr U}(F_k,V)$ by
\begin{equation}\label{Eq:augmentation-Tor-complex}
 C_\bullet(F_k,V):=C_\bullet(\mathfrak u,F_k^-\otimes V),
 \qquad X(a\otimes v)=-aX\otimes v+a\otimes Xv.
\end{equation}
In degree zero, the map $v\mapsto[1\otimes v]$ induces the algebraic identification
\begin{equation}\label{Eq:augmentation-Tor-zero}
 H_0(C_\bullet(F_k,V))
 =F_k\otimes_{\mathscr U}V\simeq V/I_{\mathfrak u}^kV.
\end{equation}
For topological $V$, the homology of this complex carries its quotient topology.

For Ext, regard $F_k$ as a left regular quotient and give its dual $F_k^\vee$ the contragredient $L$-action and the $\mathfrak u$-action $(X\lambda)(a)=-\lambda(Xa)$.
Then
\[
 \operatorname{Hom}_{\mathbb C}(F_k,V)=F_k^\vee\otimes V,
 \qquad (Xf)(a)=Xf(a)-f(Xa).
\]
The standard resolution computes $\operatorname{Ext}_{\mathscr U}^\bullet(F_k,V)$ by $C^\bullet(\mathfrak u,F_k^\vee\otimes V)$.
On the smooth side we use its cohomology quotient topologies.

\begin{proof}[Proof of Corollary~\ref{C:augmentation-comparison}]
For Tor, apply Proposition~\ref{P:finite-coefficient-comparison} with $F=F_k^-$.
Naturality gives compatibility from level $k+1$ to $k$.
For Ext, apply the cohomological assertion of the same proposition to $F_k^\vee$.
The coefficient maps are the dual injections $F_k^\vee\hookrightarrow F_{k+1}^\vee$, giving compatibility from level $k$ to $k+1$.
\end{proof}

\begin{proof}[Proof of Corollary~\ref{C:Koszul-component-closedness}]
Take $F=\bigwedge^j\mathfrak u$ or $F=(\bigwedge^j\mathfrak u)^*$.
Proposition~\ref{P:finite-coefficient-comparison} applied to $F_k^-\otimes F$ shows that $H_0(\mathfrak u,F_k^-\otimes F\otimes\pi)$ is Hausdorff.
By \eqref{Eq:augmentation-Tor-zero}, the continuous map
\[
 F\otimes\pi\longrightarrow H_0(\mathfrak u,F_k^-\otimes F\otimes\pi),
 \qquad m\longmapsto[1\otimes m],
\]
has kernel $\mathfrak u^k(F\otimes\pi)$, which is therefore closed.
\end{proof}

\subsection{Smooth Bruhat pieces}\label{S:Bruhat-piece-closedness}

Retain $J$, $Z_j$, and $J_j$ from Corollary~\ref{C:Bruhat-piece-closedness}.
The spaces $J_j$ are closed $(\mathfrak q,L)$-submodules of $J$.
The annihilator of $J_j$ in $J'$ is the support space $I_j$ of Section~\ref{S:filtrations}, so
\[
 J_j'=J'/I_j,\qquad
 (J_{j-1}/J_j)'=I_j/I_{j-1}.
\]

\begin{proof}[Proof of Corollary~\ref{C:Bruhat-piece-closedness}]
Put $V=\mathrm{CJ}_U(J_K)$ and $\mathcal J=\mathrm{CJ}^{\infty}_U(J)$.
Let $q_J:J\to\mathcal J$ be the canonical map, and put $\mathcal J_j=\overline{q_J(J_j)}$, a closed $(\mathfrak q,L)$-submodule of $\mathcal J$.
Corollary~\ref{C:CJ-comparison} identifies $\mathcal J$ with $\widehat V$.
Put $V_j=\mathcal J_j\cap V_{\mathrm{fin}}$.

Lemma~\ref{L:closed-weight-subquotients} shows that $V_j$ and $V_{j-1}/V_j$ are finitely generated over $\mathcal U(\mathfrak q)$, have Harish--Chandra weight spaces, and satisfy
\[
 \mathcal J_j\simeq\prod_\alpha((V_j)_\alpha)^\infty,
 \qquad
 \mathcal J_{j-1}/\mathcal J_j
 \simeq\prod_\alpha((V_{j-1}/V_j)_\alpha)^\infty.
\]

Applying Remark~\ref{R:weightwise-globalization} to $M=V_j$ and $M=V_{j-1}/V_j$ gives closed differential images and Casselman--Wallach homology for $\mathcal J_j$ and $\mathcal J_{j-1}/\mathcal J_j$.

To transfer these properties to the original spaces, we claim that the canonical maps induced by $q_J$ give topological $L$-isomorphisms
\begin{equation}\label{Eq:Bruhat-piece-homology-comparison}
 \begin{aligned}
 H_i(\mathfrak u,J_j)\xrightarrow{\ \sim\ }H_i(\mathfrak u,\mathcal J_j),\ 
 H_i(\mathfrak u,J_{j-1}/J_j)\xrightarrow{\ \sim\ }
 H_i(\mathfrak u,\mathcal J_{j-1}/\mathcal J_j)
 \end{aligned}
 \qquad(i\geq0).
\end{equation}

Write $\Gamma D=D^{[\mathfrak u]}$.
Under $\mathcal J'=\Gamma J'$ from \eqref{dual iso}, the annihilator of $\mathcal J_j$ is $\Gamma I_j$.
The generalized-invariant exact sequences established in the proof of Theorem~\ref{T:relative-cohomological-comparison} therefore give natural identifications of modules
\[
 \mathcal J_j'\simeq\Gamma(J_j'),
 \qquad
 (\mathcal J_{j-1}/\mathcal J_j)'
 \simeq\Gamma\bigl((J_{j-1}/J_j)'\bigr).
\]
The transposes of the canonical maps are consequently the inclusions $\Gamma D\hookrightarrow D$ for $D=J'/I_j$ and $D=I_j/I_{j-1}$.
Theorem~\ref{T:relative-cohomological-comparison} shows that they induce isomorphisms on $\mathfrak u$-cohomology.

Since the target homology is Hausdorff, Lemma~\ref{L:Frechet-duality} gives \eqref{Eq:Bruhat-piece-homology-comparison} and closed differential images for the source complexes.
Their homology is therefore Casselman--Wallach as well.

To prove that $\mathfrak u^kM$ is closed in $M$, apply the smooth assertion of Proposition~\ref{P:finite-coefficient-comparison} to $M$ and $F_k^-$.
The resulting space $H_0(\mathfrak u,F_k^-\otimes M)$ is Casselman--Wallach, hence Hausdorff.
The continuous map $M\to H_0(\mathfrak u,F_k^-\otimes M)$, $m\mapsto[1\otimes m]$, has kernel $\mathfrak u^kM$ by \eqref{Eq:augmentation-Tor-zero}.
Thus $\mathfrak u^kM$ is closed.
\end{proof}

\end{document}